\documentclass[12pt]{amsart}
\usepackage{amsmath, amscd,amssymb}
\usepackage{amsfonts}
\usepackage{braket}
\usepackage{colonequals}
\usepackage{xstring}
\input xy   %per diagrammi
\xyoption{all} %per diagrammi
\begin{document}
\numberwithin{equation}{section}

\newcommand*{\Lie}[1]{\ensuremath{\mathfrak{#1}}}
\newcommand*{\LieL}{\Lie{l}}
\newcommand*{\s}{\Lie{s}}
\newcommand*{\n}{\Lie{n}}
\newcommand*{\g}{\Lie{g}}
\newcommand*{\LieH}{\Lie{h}}
\newcommand*{\LieG}{\Lie{g}}
\newcommand*{\LieK}{\Lie{k}}
\newcommand*{\LieN}{\Lie{n}}
\newcommand*{\LieM}{\Lie{m}}
\newcommand*{\LieV}{\Lie{v}}
\newcommand*{\LieSL}[1]{\ensuremath{\mathfrak{sl}\!\left(#1\right)}}
\newcommand*{\SL}[1]{\ensuremath{SL\!\left(#1\right)}}
\newcommand*{\uSL}[1]{\ensuremath{\widetilde{SL}\!\left(#1\right)}}
\newcommand*{\de}{\partial}
\newcommand*{\R}{\mathbb R}
\newcommand*{\of}[1]{\ensuremath{\left(#1\right)}}
\newcommand*{\Aut}[1]{\ensuremath{\operatorname{Aut}\!\of{#1}}}
\newcommand*{\Fix}[1]{\ensuremath{\operatorname{Fix}\!\of{#1}}}
\newcommand{\C}{\mathbb C}
\newcommand{\B}{\mathbb B}

\newcommand*{\defeq}{\mathrel{\vcenter{\baselineskip0.5ex \lineskiplimit0pt
                     \hbox{\scriptsize.}\hbox{\scriptsize.}}}%
                     =}

\newcommand*{\LieDer}{\ensuremath{\EuScript L}}

\newtheorem{theorem}{Theorem}[section]
\newtheorem{lemma}[theorem]{Lemma}
\newtheorem{proposition}[theorem]{Proposition}
\newtheorem{corollary}[theorem]{Corollary}

\theoremstyle{definition}
\newtheorem{definition}[theorem]{Definition}
\newtheorem{example}[theorem]{Example}

\theoremstyle{remark}
\newtheorem{remark}[theorem]{Remark}
\numberwithin{equation}{section}

%tentative title
\title[Homogeneous holomorphic foliations]{Invariant holomorphic foliations on Kobayashi hyperbolic homogeneous manifolds}
\author[F. Bracci]{Filippo Bracci$^\dag$}
\address{F. Bracci \and A. Iannuzzi: Dipartimento Di Matematica\\
Universit\`{a} di Roma \textquotedblleft Tor Vergata\textquotedblright\ \\
Via Della Ricerca Scientifica 1, 00133 \\
Roma, Italy} \email{fbracci@mat.uniroma2.it, iannuzzi@mat.uniroma2.it}
\author[A. Iannuzzi]{Andrea Iannuzzi}
\author[B. McKay]{Benjamin McKay}
\address{B. McKay: University College Cork\\ National University of Ireland\\ Cork, Ireland} \email{b.mckay@ucc.ie}
\date\today
\thanks{$^{\dag}$Supported by the ERC grant ``HEVO - Holomorphic Evolution Equations'' n. 277691.}

\begin{abstract}
Let \(M\) be a  Kobayashi hyperbolic
homogenous manifold.
Let \(\mathcal F\) be a holomorphic foliation on \(M\) invariant under a transitive group \(G\) of biholomorphisms.
We prove that the leaves of \(\mathcal F\) are the fibers of a holomorphic \(G\)-equivariant submersion \(\pi \colon M \to N\) onto a \(G\)-homogeneous complex manifold \(N\).
We also show that if $\mathcal Q$ is an automorphism family of  
 a hyperbolic convex (possibly unbounded) domain  \(D\) in $\C^n$, then the
fixed point set of \(\mathcal Q\) is
either  empty or a connected complex submanifold of \(D\).
\end{abstract}



\maketitle
\tableofcontents

%----------------------------------------------------------
%------------------------INTRO-----------------------
%----------------------------------------------------------

\section{Introduction}

The existence of a nontrivial holomorphic foliation on a complex
manifold \(M\)
gives rise to restrictions on the geometry of the manifold itself. For instance, by the Baum--Bott index theorem \cite{BaBo}, the existence of such a foliation 
implies the vanishing of certain characteristic classes of \(M\).
As a consequence, complex projective spaces do not admit
nonsingular holomorphic foliations. Recently,
M.~Brunella,  M.~McQuillan and L.~G.~Mendes (see \cite{Bru})  gave  a birational classification of (singular) holomorphic foliations on projective surfaces in the spirit of the Enriques--Kodaira classification.
On noncompact manifolds all characteristic classes tend to vanish and a
similar result is not to be expected without further assumptions.

Let \(G\) be a Lie group acting by biholomorphisms on \(M\) preserving a holomorphic foliation \(\mathcal F\), i.e. if \(F\) is a leaf of \(\mathcal F\), then so is \(g \cdot F\), for all \(g\in G\).
 In this context, A.~Behague and B.~Sc\'ardua \cite{BeSc} gave a complex version of a
classical result of D. Tischler \cite{Tis}. 
Namely, they showed that a
holomorphic foliation with closed leaves
 which is  invariant under a holomorphic transverse action of a complex Lie group of
 dimension equal to the codimension of the foliation,
 is given by a holomorphic fibration.

In the case of a transitive action of a  Lie group \(G\), the leaves
of a \(G\)-invariant foliation are nonsingular but might not be closed.
For instance a
foliation defined by a generic translation vector field on a complex
torus \(\,\C^n/\Lambda\), or
on  \( (\C^*)^n\), has no closed leaves.

Here we assume the manifold to be homogeneous and Kobayashi hyperbolic which, by a result of K. Nakajima \cite{Nak}, implies that $M$ is biholomorphic to a homogenous Siegel domain of type II.
The aim of this paper is to prove the following theorem.%----------------------------------------------------------

\begin{theorem}\label{mainr}
Let  \(M\) be a hyperbolic homogeneous complex manifold and $G$ a group of 
automorphisms acting transitively on $M$. Let \(\mathcal F\) be a $G$-invariant holomorphic foliation on \(M\). Then there exists a \(G\)-homogeneous complex manifold \(N\) and a \(G\)-equivariant holomorphic submersion \(\pi \colon M \to N\) such that the leaves of \(\mathcal F\) are the fibers of \(\pi\).
\end{theorem}

%----------------------------------------------------------

The base $N$ might not be hyperbolic
(cf. Example \ref{BALL}). 
The hyperbolicity of \(M\) implies that 
the automorphism group \(\Aut M\) of \(M\) is a Lie group acting properly on \(M\).
The foliation \(\mathcal F\) is  invariant with respect to the closure of \(G\) in 
\(\Aut M \). Thus,  without loss of generality  we may assume  \(G\) to be closed in \({\rm Aut}(M)\),
implying that $G$ is a Lie group acting properly on $M$.

Examples of equivariant submersions whose base and total space are
hyperbolic homogeneous Siegel domains  are given,
{\sl e.g.} in \cite{Mie}, Prop. 5.6,  where C. Miebach remarked that 
such submersions do not need to be holomorphic fiber bundles
(cf. Example \ref{ISHI}).  In fact, if the submersion $\pi$ in the above theorem is a holomorphic fiber bundle then \(M\) is biholomorphic to a product  \(N \times F\) of hyperbolic, homogeneous Siegel domains
and \(G\) is a subgroup of \(\Aut N \times \Aut F\) (cf. Lemma \ref{trivial-fib}).

If the foliation is not invariant, then hyperbolicity is not  a sufficient condition
for the leaves to be closed. In Section 4
we  present an example (Example \ref{forn})  suggested by John Erik Forn\ae ss of a nonsingular holomorphic foliation on the unit ball \(\B^2 \subset \C^2\) having some nonclosed leaves.

The paper is organized as follows. In Section \ref{BF} we introduce
``uniform Bochner-Frobenius local coordinates'' and present some preliminary 
material. In Section \ref{LEAVES} we  show that the leaves of \(\mathcal F\) are closed. In Appendix \ref{Lie} we give a different proof of this fact by exploiting the
 Lie group structure of the automorphism group of a Siegel domain.
 In Section~\ref{Trivial} we prove Theorem \ref{mainr} and discuss
some examples.

In Section \ref{A-fixed} we consider the foliation $\mathcal F^K$ 
induced on  the fixed point set $M^K$ of the isotropy subgroup $K$ of $G$ at
one point. The set  
$M^K$ is a  connected homogenous complex  submanifold of $M$. 
We show that every leaf of $\mathcal F^K$ is the intersection of a leaf of $\mathcal F$ with $M^K$, i.e. such an intersection is connected. The result is based on the following proposition, which might be of interest on its own:

\begin{proposition}
\label{lemma:connected}
Let \(D\) be a hyperbolic convex $(\,$possibly unbounded$\,)$ domain in \(\C^n\)
and let  \( \mathcal Q\) be a family of  automorphisms of \(D\). 
Then the fixed point set \(D^ {\mathcal Q}\)
 of \(\mathcal Q\) is either empty
or a connected complex submanifold of \(D\).
\end{proposition}

Note that if  \(D^ {\mathcal Q}\)
  is not empty, then \(\mathcal Q\) is relatively compact
in $\Aut D$.
Using Bochner coordinates one easily sees that each connected component of 
\(D^ {\mathcal Q}\) is  a complex submanifold of the domain. The main issue in the above proposition
is the connectness of  \( D^{\mathcal Q}\).

%----------------------------------------------------------

\medskip
\noindent
{\bf Acknowledgments} We wish to thank 
John Erik Forn\ae ss for suggesting
Example \ref{forn}. We are also grateful to the referee for her/his accurate comments
and  for pointing out a gap in a previous version of the paper.

%----------------------------------------------------------
%-----------------------BF----------------
%----------------------------------------------------------

\section{Preliminaries}
\label{BF}

Let \(M\) be a hyperbolic  complex manifold. We refer to \cite{Kob} for the definition of Kobayashi distance, hyperbolic manifolds and their properties. Let \(G\)  be a  
connected, closed subgroup of \({\rm Aut}(M)\). Then  the \(G\)-action is proper
(cf. \cite{Kob}, Thm. 5.4.2). In particular 
every isotropy subgroup  
\[
G_p \defeq \Set{g\in G \ : \ g\cdot p=p}
\]
is compact. 
Denote by \(\B^m\) the unit radius Euclidean ball in \(\C^m\) centered at \(0\).

%----------------------------------------------------------

\begin{lemma}
\label{BF-coordinates}
Let \(M \cong G/K\) be a 
hyperbolic, homogeneous $n$-dimensional complex manifold and let  \(\mathcal F\) be a  $k$-dimensional, \(G\)-invariant  holomorphic foliation on \(M\). Then for any \(p\in M\) there exist holomorphic local coordinates \((U_p, \psi_p)\), with $U_p$ a $G_p$-invariant neighbourhood of \(p\), such that:
  \item{\rm (i)}  \(\psi_p(U_p)=\B^{n-k}\times \B^k\),
  \item{\rm (ii)} the plaques of  \(\mathcal F\) in \(U_p\) 
correspond to   $\{\,(z,w) \in \B^{n-k}\times \B^k \ : \ z= \textrm{const}\, \}$,
 \item{\rm (iii)}  there is a faithful representation $L:G_p \to U(n-k) \times U(k)$
 so that  for every $g \in G_p$  and $q \in U_p$ one has 
  $\psi_p(g \cdot q)= L(g)\cdot \psi_p(q)$. \end{lemma}

\begin{proof}
Since the foliation is nonsingular near \(p\), there are local holomorphic coordinates \((z,w)\) centered at \(p\) (called ``Frobenius coordinates'') 
in which the plaques of \(\mathcal{F}\) are \(\{z=\textrm{const}\}\).
Since \(G_p\) is compact, averaging any inner product on $T_pM$ over Haar measure
gives a \(G\)-invariant Riemannian metric on \(M\) inducing the standard topology.
Since \(G_p\) fixes \(p\), we can pick a \(G_p\)-invariant open set $U$ (for example, a ball in such a metric) sufficiently small as to lie in the domain of the coordinates.
So without loss of generality, we can assume that  \(U \subset \C^n\) is a bounded domain, \(p=0\) and the foliation \(\mathcal{F}\) on $U$ is given by \(\{z=\textrm{const}\}\).

Identify \(\C^n\) with its tangent space at the origin and for \(g\in G_p\) consider the linear operator on \(\C^n\)
defined by \(L(g)=dg_p\). Then  \(L(G_p)\) is a compact group
of linear  transformations of \(\C^n\) permuting the affine subspaces  \(\{z=\textrm{const}\}\).
By choosing an \(L(G_p)\)-invariant Hermitian inner product on \(\C^{n}\) we realize
\(L(G_p)\) as a closed subgroup of  the unitary group \(U(n)\).
 Moreover, since the subspace  \(\{z=0\}\) is invariant, one has
 \(L(G_p)\subseteq U(n-k)\times U(k)\).

We conclude by showing that such an \(L(G_p)\)-action on \(\C^n\) is a local linearization of the \(G_p\)-action on \(M\).
Following  \cite{Boc}, p. 375, consider the Haar measure  \(\mu\)
on \(G_p\) and  define a holomorphic map on \(U\)  by
\[
\psi_p(z,w) \defeq
\int_{G_p} L(h)^{-1} \cdot(
 h \cdot (z,w) )d\mu(h).
\]
By construction \((d\psi_p)_p=\operatorname{Id}\) and,
shrinking \(U\) if necessary, \(\psi_p\) is a biholomorphism onto its image.
Also, \(\psi_p \circ g= L(g) \circ  \psi_p\) for all \(g\in G_p\).
Note that
\(\psi_p\) preserves \(\mathcal{F}\), since both  \(g\) and \(L(g)\) permute the affine subspaces \(\{z=\textrm{const}\}\). 
Finally,  arrange so that the image of  \(\psi_p\) coincides with a copy of \(\B^{n-k}\times \B^k\).
\end{proof}

The above coordinates, which we refer to as 
{\sl Bochner--Frobenius local coordinates}
are {\it uniform}  with respect to
the Kobayashi distance of \(M\) in the following sense.

%----------------------------------------------------------

\begin{lemma}\label{BFuniform}
Let  \(M\cong G/K\) be a  hyperbolic,  homogeneous complex manifold and let
\(\mathcal F\) be a  \(G\)-invariant  holomorphic foliation
on \(M\).
Then there exists \(r>0\) such that every point \(p \in M\) lies in a system of
Bochner--Frobenius local coordinates, centered at \(p\), and defined on the Kobayashi ball of radius \(r>0\) centered at \(p\).
\end{lemma}

\begin{proof}
Since \(M\) is complete hyperbolic (see \cite{Kob}, Thm. 3.6.22), the topology of \(M\) coincides with that defined by the Kobayashi distance.
Thus,  given \(p_0 \in M\)
there exists \(r > 0\) such that the
Kobayashi ball centered at \(p_0\) of radius \(r\) lies
inside the domain of Bochner--Frobenius coordinates $(U_{p_0}, \psi_{p_0})$.
Let \(p \in M\) and \(g\in G\) be such that \(g \cdot
p_0=p\). 
It is easy to check that \((g \cdot U_{p_0},\psi_{p_0} \circ g^{-1})\) are Bochner--Frobenius
local coordinates centered at \(p\). Then the result follows by
recalling that \(g\) is an isometry for the Kobayashi distance.
\end{proof}

%----------------------------------
%Let  \(M\) be a  hyperbolic,  homogeneous complex manifold and let
%\(\mathcal F\) be a  holomorphic foliation
%on \(M\) invariant under the action of a group \(G\).  
%An \emph{invariant local transverse} to \(\mathcal F\) at a point \(p \in M\) is a \(G_p\)-invariant %embedded \(k\)-codimensional complex submanifold \(S_p \subset M\) containing \(p\) so that %each leaf of \(\mathcal F\) either does not intersect \(S_p\) or intersects \(S_p\) transversally.
%By Lemma \ref{BF-coordinates}, the set
 %\(
%S_p = \psi_p^{-1}( \B^{n-k}\times \{0\} )
%\)
%is an invariant local transverse at \(p\).

%----------------------------------------------------------

%\begin{lemma}\label{BF-transverse}
%Let  \(M\cong G/K\) be a  hyperbolic,  homogeneous complex manifold and let
%\(\mathcal F\) be a  \(G\)-invariant  holomorphic foliation
%on \(M\).
%Let \(F_p \subset M\) be the leaf of \(\mathcal F\) through a point \(p \in M\) and \(S_p \subset M\) be an  invariant local transverse through \(p\).
%If $G_p$ is connected, then \(F_p \cap S_p\) is contained in the fixed point set of  \(G_p\).
%\end{lemma}

%\begin{proof}
%Let \(q\in F_p\cap S_p\). Since both submanifolds are \(G_p\)-invariant, 
%the \(G_p\)-orbit of \(q\) is also contained in \(F_p\cap S_p\).
%Moreover,  \(F_p\) and \(S_p\) intersect transversally, implying that
%\(G_p \cdot q\)  cannot be positive dimensional. 
%Finally, if \(G_p\) is connected then such an orbit is connected,
%therefore it consists of a fixed point.
%\end{proof}

The following fact will be used in Section \ref{Trivial} in
order to show that the fibration onto the leaf space of the foliation is 
a holomorphic submersion.

%----------------------------------------------------------

\begin{proposition}
\label{ALMOSTCPLX}
Let \(M\) be a complex manifold and $\mathcal F$ a holomorphic foliation on $M$. Assume there exists a real manifold \(N\) and a $C^\infty$ submersion \(\pi \colon M \to N\)  such that the leaves of $\mathcal F$ are the fibers of $\pi$.
Then there exists a unique complex structure on \(N\) so that \(\pi\) is holomorphic.
\end{proposition}

\begin{proof}
We say that a complex-valued function on an open subset of \(N\) is an \emph{up-function} if its composition with \(\pi\) is holomorphic.
A complex structure on \(N\) is \emph{down} if the  functions which are holomorphic with respect to  that complex structure are precisely the up-functions, on every open subset of \(N\).

Note that any two complex structures are equal if they have the same holomorphic functions on every open set. Therefore a down complex structure on \(N\) is unique, if it exists. The problem of existence is local: should we prove existence of a down complex structure on each open set in a covering of \(N\), then by uniqueness these complex structures agree where any two are defined, so they glue together to a unique down complex structure on \(N\).

Fix a point \(p_0 \in M\). Choose local Frobenius coordinates $(U, (z,w))$ centered at $p_0$,
so that the plaques of $\mathcal F$ are given by $\{z=\hbox{const}\}$. The complex submanifold
$Z:=\{(z,w) :w=0\}$ of $M$ is transverse to the foliation and, 
by shrinking $U$ if necessary, the map \(\pi\) restricts to \(Z\) to be a local diffeomorphism \(\left.\pi\right|_Z \colon Z \to \pi(Z)\).
This diffeomorphism defines  a (integrable) complex structure on the open set $\pi(Z)\subset N$.
In this complex structure, a holomorphic function on \(\pi(Z)\) is precisely one 
whose pull-back is holomorphic on \(Z\).
In particular, up-functions on open subsets of \(\pi(Z)\) are holomorphic.

In order to conclude the proof, we need to show that every function  on an open subset $W$
of \(\pi(Z)\), which is holomorphic with respect to such a complex structure, is in fact an up-function. 
The pull-back \(f:\pi^{-1}(W) \to \C \) of such a function is constant along the fibers of $\pi$ over $W$ and holomorphic on $\pi^{-1}(W) \cap Z$. With respect to the above Frobenius coordinates, 
 $f$ is given by $(z,w)\mapsto f(z,w)=f(z,0)$. Thus it is  holomorphic on $\pi^{-1}(W) \cap U$. 

In order to check  that $f$ is holomorphic on $\pi^{-1}(W)$, let $q_1\in \pi^{-1}(W)$. Then there exists a unique point $q_0\in Z$ such that $q_1\in \pi^{-1}(\pi(q_0))$. By assumption,  the fiber $\pi^{-1}(\pi(q_0))$ is the leaf  of $\mathcal F$ through $q_0$, hence it is connected. Therefore there exists a continuous curve \(\rho \colon [0,1] \to \pi^{-1}(\pi(q_0))\)  such that $\rho(0)=q_0$ and $\rho(1)=q_1$. 
The set
\[
A:=\{t\in [0,1] : f \hbox{ is holomorphic in a neighborhood of\ }\rho(t) \ {\rm in}\  M\}
\]
is clearly open in $[0,1]$. Moreover $0 \in A$, since $f$ is holomorphic on $U \ni \rho(0)$. 
Suppose that \(A \ne [0,1]\); take \(t_1 \in [0,1]\) to be the smallest real number not in \(A\).
Take a new Frobenius coordinate system $(U_1,(z,w))$ centered at  \(\rho(t_1)\) such that 
$U_1 \subset \pi^{-1}(W)$.
Then, for $t_2$ close enough to \(t_1\),  with $t_2<t_1$, the point $\rho(t_2)$ belongs  to 
$U_1$. Thus $\rho(t_2) = (0,w_2)$, for some $w_2$ in
$\C^{\dim\mathcal F}$.

By the definition of $t_1$, the function $f$ is holomorphic near $(0,w_2)$.
Since in the Frobenius coordinate system $(U_1,(z,w))$ the function
 \(f\) is independent of \(w\), it follows that \(f(z,0)=  f(z,w_2)\).
This implies that $f(z,w)= f(z,0)$
is holomorphic in a  neighborhood of  \(\rho(t_1)=(0,0)\),
giving a contradiction. Thus $A=[0,1]$ and \(f\) 
is holomorphic in a neighborood of $q_1$ for every $q_1 \in \pi(W)$, i.e.
$f$ is holomorphic on $\pi^{-1}(W)$.
Hence the defined complex structure on $\pi(Z)$ is down, concluding the proof.
\end{proof}

\section{Leaves are closed}
\label{LEAVES}
In this section we show that under the assumptions of Theorem \ref{mainr}
the leaves of the foliation are closed.

\begin{proposition}
\label{CLOSED}
Let \(G\) be a Lie group acting transitively on 
a hyperbolic, complex manifold \(M\)
and let \(\mathcal F\) be a  \(G\)-invariant holomorphic foliation  on
 \(M\). Then the leaves of \(\mathcal F\) are closed.
\end{proposition}
\begin{proof}

Let $p$ be a point in $M$ and
assume by contradiction that the leaf $F_p \in \mathcal F$ through $p$ is not closed. 

\smallskip
\noindent
{\it Claim.} There exists Bochner-Frobenius coordinates $(U_p,\psi_p)$ at $p$ 
and a sequence $\{p_n\} \subset U_p \cap F_p$ such that $p_n \to p$ and,
for every $n \in \mathbb N$, the points $p_n$
do not belong to the plaque of $F_p$ in $U_p$ containing $p$.

\smallskip
\noindent
{\it Proof of claim.} 
 The plaques of $\mathcal F$ on $U_p$ are given by $\{(z,w) \in \B^{n-k}
 \times \B^k\ :\ z=\hbox{const}\}$, and $p$ belongs to the plaque $\{z=0\}$.  If a sequence as in the claim does not 
 exist then, by shrinking $U_p$ if necessary, we may assume that $F_p\cap U_p$
consists of a  unique plaque, namely $\{z=0\}$. Since $F_p$ is not closed, 
 there exists a sequence of points $q_n\in F_p$ such that $q_n\to q$, with $q\not\in F_p$.
By homogeneity (see Lemma \ref{BFuniform}), there 
 exists Bochner-Frobenius coordinates $(U_{q_n},\psi_{q_n})$ at $q_n$ such that
 $F_p\cap U_{q_n}$ consists of a  unique plaque and, for $n$ large enough,
  $q \in U_{q_n}$. In particular  $F_p\cap U_{q_n}$ is closed in $U_{q_n}$
  and since  $q_m \in F_p \cap U_{q_n}$, for $m$ large enough,
 this gives a contradiction and proves the claim.
   \smallskip
 
The leaf $F_p$ can be regarded as an immersed complex submanifold
$\varphi:N\to M$, with  $\varphi(N)=F_p$.
Let $p_n$ be as in the claim, $\zeta, \,\zeta_n\in N$ such that  $\varphi(\zeta)=p$ and $\varphi(\zeta_n)=p_n$. Since $p_n$ does not belong to the same plaque of $p$ and $\varphi$ is an immersion, it follows that $\zeta$ is not in the closure of the set 
$\{\zeta_n\}$. Pick $g_n\in G$ such that $g_n \cdot p=p_n$ and, 
consequently, $g_n \cdot F_p=F_p$.
Note that $g_n$ maps plaques of $\mathcal F$ on a neighborhood $U$
onto plaques of $\mathcal F$ on $g_n \cdot U$, and the plaques are open in $F_p \cong N$. Thus, $h_n:=\varphi^{-1}\circ g_n \circ \varphi:N\to N$ is a biholomorphism  of $N$ such that $h_n \cdot \zeta=\zeta_n$. 

Since $M$ is complete hyperbolic and hence taut (see, {\sl e.g.} \cite{Kob}) and $p_n=g_n\cdot p \to p$,  we can assume with no loss of generality that $\{g_n\}$ converges uniformly on compacta to some $g\in G_p$. Hence $h_n \to  h:=\varphi^{-1}\circ g \circ \varphi:N\to N$, with $h$ an automorphism of $N$ such that $h\cdot \zeta=\zeta$. This implies 
 $\zeta_n =h_n \cdot \zeta \to \zeta$ and gives a contradiction.
\end{proof}

%----------------------------------------------------------
We conclude this section with an example of a nonsingular holomorphic foliation on the unit ball \(\mathbb B^2\) of
\(\C^2\) having nonclosed generic leaves.

\begin{example}\label{forn}
Let \(X\) be the real vector field in \(\R^2\) given by
\[
X(x,y)\defeq -y\frac{\de}{\de x}+x\frac{\de}{\de y}
+\left(\frac{1}{2}-x^2-y^2\right)
\left(x\frac{\de}{\de x}+y\frac{\de}{\de y}\right).
\]
Then the circle \(C:= \{x^2+y^2=\frac{1}{2}\}\) is a closed integral curve of \(X\).
All nearby leaves accumulate to that curve. Indeed for each \((x,y)\) close enough to \(x^2+y^2=\frac{1}{2}\) the standard
scalar product
\[
\left< X(x,y),x\frac{\de}{\de x}+y\frac{\de}{\de y}\right>
\]
is positive if \((x,y)\) belongs to the disc of radius \(1/2\), and negative otherwise.
Hence \(X\) pushes towards $C$.

Let \(Z\) be the complexification of \(X\) in the \((z,w)\)-plane, {\sl i.e.}
 \(Z(z,w)\) is given by replacing \(x\) with \(z\) and \(y\) with \(w\). Let \(\mathcal F\) be the holomorphic foliation defined by \(Z\). A direct computation shows that the only singularity of \(\mathcal F\) is \((0,0)\). Moreover, the complex conic
 \(C^\C:=\Set{(z,w)\in \C^2|z^2+w^2=1/2}\) is a closed leaf of \(\mathcal F\) and nearby leaves accumulate to it.
Consider the compact,  polynomially convex set
\[
K \defeq
\Set{(z,w)\in \C^2| z^2+w^2=1/2, |z|^2+|w|^2\leq 2}
\]
not containing the origin. By \cite{For}, Prop. 2.1, given a small positive \(\varepsilon > 0\),
 there exists an automorphism \(\Phi\in {\rm Aut}(\C^2)\) such that \(\Phi(0,0)=(2,0)\) and \(|\Phi(z)-z|<\varepsilon\) for all \(z\in K\).

Consider the foliation \(\tilde{\mathcal F}\) given by the restriction to \(\B^2\) of the holomorphic foliation \((\Phi^{-1})^\ast \mathcal F\). By construction the foliation
\(\tilde{\mathcal F}\) has no singularities in \(\B^2\) and,
for $\varepsilon$ small enough, $ \Phi(C)$ is contained 
in $\B^2$.
Then the connected component of \(\Phi(C^\C)\cap \B^2\) 
containing $ \Phi(C)$  is a closed leaf of 
\(\tilde{\mathcal F}\) and nearby leaves accumulate on it, thus they are not closed.
\end{example}

%----------------------------------------------------------
%---------------Proof-----------------
%----------------------------------------------------------

\section{The equivariant submersion}
\label{Trivial}

\begin{proof}[Proof of Theorem \ref{mainr}]
Without loss of generality we may assume that \(G\) is a connected, closed subgroup of \(\Aut M\).
The leaf \(F\) through the base point is closed by Proposition \ref{CLOSED}, therefore so is its stabilizer \(L\) in \(G\). 
Let \(N \defeq G/L\), a homogeneous manifold with \(G\)-equivariant fiber bundle map \(\pi \colon M=G/K \to N=G/L\).
By Proposition~\ref{ALMOSTCPLX}, there is a unique complex structure on the leaf space \(N\) such that \(\pi\) is holomorphic. Since the real analytic $G$-action on $M$ is by biholomorphisms,
so is  the $G$-action on $N$.
\end{proof}

%----------------------------------

Let \(M\) be a hyperbolic \(G\)-homogeneous manifold, let \(N\) be a \(G\)-homogeneous manifold and let
\(\pi \colon M \to N\) be a \(G\)-equivariant holomorphic submersion. The \(G\)-equivariance implies that  such a submersion 
is a smooth fiber bundle. However,  in general \(\pi\)  does not admit any local 
holomorphic trivialization, as noticed in the case of certain equivariant
submersions by C. Miebach in \cite{Mie}, Prop.5.6 and Rem. p. 347.
The fibers of Miebach's submersions are biholomorphic to
balls and the bases of his submersions are hyperbolic.
In this section we construct an example where the base of the submersion is 
not hyperbolic. We first collect some useful facts.

%---------------------Lemma

\begin{lemma}\label{trivial-fib}
Assume that the \(G\)-equivariant submersion \(\pi \colon M \to N\) in 
Theorem \ref{mainr} is a holomorphic fiber bundle.
Then \(\pi\) is a trivial holomorphic fiber bundle,
\(M\) is biholomorphic to a product  \(N\times F\)
of  hyperbolic, homogeneous Siegel domains
and \(G\) is a subgroup of \(\Aut N \times \Aut F\).
\end{lemma}

\begin{proof}
Since the total space of the holomorphic fiber bundle \(\pi \colon M \to N\)
is  hyperbolic, the base
\(N\)  is also hyperbolic by a result of  S. Nag \cite{Nag}.
As a consequence, \(N\) is biholomorphic to a homogenous
Siegel domain of type II \cite{Nak}. 
In particular it is
simply connected. Then, Royden 
 \cite{Roy}, Cor. 1, implies that
\(\pi\) is a trivial holomorphic fiber bundle and 
consequently \(M\) is biholomorphic to \(N\times F\).

We are left to show that \(G\) is contained in the product 
\(\Aut N \times \Aut F\) (cf. \cite{Kob}, Cor. 5.4.12, for 
 $G$  connected).
 
Write \(M=N\times F\) and for \(g\) in \(G\) let 
$\,g \cdot (z,w)= (\alpha_g(z,w), \beta_g(z,w))\,.$
Since \(g\) preserves the leaves \(\{*\} \times F\)
of \(\mathcal F\), it follows that \(\alpha_g\) does not depend on \(w\). 
Then one checks that the map \(G \to \Aut N\), given by \(g \to \alpha_g\) is
a group homomorphism.

Fix \(z\in N\) and consider the map
\(\beta_g (z,*) \colon F \to F\), defined by \(w \to \beta_g (z,w)\).
Let \(g^{-1} \cdot (z,w)=(\alpha_g^{-1}(z),
\beta_{g^{-1}}(z,w))\).
Since \(g^{-1} g \cdot (z,w)=(z,w)\),
 it follows that 
 $$\beta_{g^{-1}}(\alpha_g(z), \beta_g(z, *))={\rm Id}_F \,.$$
Hence, \(\beta_{g^{-1}} (\alpha_g(z),*)\) is a left inverse  of \(\beta_g(z, *)\).
Similarly, one checks that it is  a right inverse as well.
Thus \(\beta_g (z,*)\) is an automorphism of \(F\). Moreover
  the map
\[N \times F \to F, \quad \quad (z,w) \to \beta_g (z,w)\,,\]
is holomorphic. Then Proposition  in \cite{Roy}
applies to show that \(\beta_g\)  does not depend
on \(z\) and one checks  that the map 
\[G \to \Aut N \times \Aut F, \quad \quad g \to (\alpha_g, \beta_g) \,,\]
is a group homomorphism with trivial kernel.
 \end{proof}

%---------------Example 

\begin{example}\label{BALL}
Let \(M=\B^2\) be the unit ball of \(\C^2\) and \(G\) the \(5\)-dimensional
isotropy subgroup of \(\Aut {\B^2} \cong SU(2,1)\) at a boundary
point \(q \in \partial \B^2\). The group  \(G\) is of the form 
\(T S\), where \(T\) is a one-dimensional torus of \(SU(1,2)\)
and \(S=AN\) is the solvable factor of an Iwasawa decomposition \(KAN\)
of \(SU(1,2)\). In particular \(G\)
acts transitively on \(\B^2\) and 
leaves invariant the foliation whose leaves are 
the complex geodesics whose closure contains 
\(q\) (cf. \cite{Aba}, Cor. 2.6.9, p. 308).
Moreover, the space of leaves is biholomorphic to the space of 
complex lines through \(q\) which are not tangent to \(\B^2\), {\sl i.e.} to \(\C\).

In order to give a  simple realization of this construction,
it is convenient to consider the hyperbolic model 
 \(\mathbb H^2 \defeq \{ \,(z,w) \in \C^2 \ : \ {\rm Im}\, w >
|z|^2\}\) of \(\B^2\)
embedded in \(\mathbb P^2\) via the map
\[ \mathbb H^2 \to \mathbb P^2, \quad (z,w) \to [z:w:1]\,.\]
Any complex geodesic of \(\mathbb H^2\) whose closure contains the boundary point \([0:1:0]\)
is given by \(\{z=const\}\). These complex geodesics define a holomorphic foliation on \(\mathbb H^2\)
which is invariant with respect to the
\(5\)-dimensional isotropy group \(G\) of \(\Aut {\mathbb H^2}\) at   \([0:1:0]\).
The isotropy of \(G\) at \((0,i)\) is given by \(\,
\{\,(z,w) \to (e^{i\theta}z,w)\ : \ \theta \in \R\,\}\) and \(G\) contains the
solvable subgroup generated by the elements of the form
$$\begin{aligned}
&(z,w) \to (z, w +t)\,, \quad \quad & &(z,w) \to (e^tz, e^{2t}w)\,,\cr
&(z,w) \to (z+t, w+2itz+ it^2)\,, \quad  \quad &&(z,w) \to (z+it, w+2tz+ it^2)\,,
\end{aligned}$$
for \(t \in \R\) (cf. \cite{Mie}, Sect. 4).
In particular, \(G\) acts transitively on \(\mathbb H^2\).
 
Every leaf of the foliation is biholomorphic to the unit
disc  in \(\C\).
The \(G\)-equivariant submersion of Theorem \ref{mainr} is given by the 
projection
\[\pi \colon \mathbb H^2 \to \C\,,  \quad\quad (z,w) \to z\,,\]
and, as a consequence of  the above lemma, it is not a holomorphic fiber bundle.
 \end{example}

 \begin{example}
 \label{ISHI}
 Let \(M\) be a hyperbolic homogenous manifold. 
 By a result of H.~Ishi \cite{Ish} the isotropy 
subgroups of  \(\Aut M\) are at least one
 dimensional. In fact, such a result is sharp. As an example,
 consider the tube domain \(M_V\) over the Vinberg cone
 \[V \defeq \{(y_1,y_2,y_3,y_4,y_5) \in \R^5 \ : \ y_3>0, \ 
 y_1y_3-y_4^2>0, \ y_2y_3-y_5^2>0\}\,.\]
Then, by \cite{Gea}, Prop. 2.2, the isotropy subgroup of \(\Aut {M_V}\) at any point of \(M_V\)
 is 1-dimensional.
Let \(G\) be a transitive subgroup of \(\Aut {M_V}\) and \(\mathcal F\) a \(G\)-invariant 
foliation (e.g. the fibration constructed in \cite{Mie}).
 Then the induced equivariant submersion  is not a holomorphic fiber bundle.
 If it were, by the above lemma the one-dimensional isotropy subgroup \(K\) of
 \(\Aut {M_V}\)  at a point \((z_0,w_0) \in N \times F \cong M\) would contain the product of the isotropy subgroups
 of  \(\Aut N\) at \(z_0\) and \(\Aut F\) at \(w_0\). Thus, by Ishi's result, \(K\)
 would be at least two-dimensional, giving a contradiction. 
  \end{example}
  
 \section{Fixed point sets of isotropy subgroups}
 \label{A-fixed}

In this section we first prove Proposition \ref{lemma:connected}.

%----------------------------------------------------------

\begin{lemma}\label{lemma:retraction}
Let \(D\) be a hyperbolic convex (possibly unbounded) domain in \(\C^n\).
The Kobayashi balls of \(D\) are bounded convex domains.
\end{lemma}

\begin{proof}
Let \(\B_N\) be the ball in \(\C^n\) with center \(0\) and radius \(N \in \mathbb N\).
Let \(D_N \defeq D\cap \B_N\). Note that \(D_N\) is a bounded convex domain for
 all \(N\in \mathbb N\). Let \(k_{D_N}\) denote the Kobayashi distance of \(D_N\). Since
 \(D\) is the increasing union of the \(D_N\), it follows that \(k_D=
 \lim_{N\to \infty} k_{D_N}\). For every \(N\), the Kobayashi distance
 \(k_{D_N}\) is a convex function (see \cite{Aba}, Prop. 2.3.46). 
Passing to the limit, \(k_D\) is a convex function as well. 
By \cite{BrSa}, a hyperbolic convex domain is complete hyperbolic.
Hence, the Kobayashi balls of \(D\) are bounded convex domains in
 \(\C^n\).
 \end{proof}

%----------------------------------------------------------

\begin{proof}[Proof of Proposition \ref{lemma:connected}]
Assume that 
\[ D^{\mathcal Q} \defeq \{z \in D \ : \ g \cdot z = z \ {\rm for \ all \ }  g   \in \mathcal Q \}\]
is nonempty and let \(z \in D^{\mathcal Q}\). Then \(\mathcal Q\) is contained inside
the  isotropy subgroup of \(z\) in \(\Aut D\) which, by the hyperbolicity of \(D\), 
is compact. Using Bochner's local coordinates  
one sees that \(D^{\mathcal Q} \) is a closed complex submanifold of \(D\) (cf. \cite{Aba}, Cor 2.5.10).
Moreover, there exists \(s(z)>0\) such that for all \(0<r<s(z)\) the intersection of
 the Kobayashi ball \(B(z,r)\) with  \(D^{\mathcal Q}\) lies in 
 the connected component \((D^{\mathcal Q})^z\) of 
\(z\) in \({D^{\mathcal Q} }\). 

Assume by contradiction that \(D^{\mathcal Q}\)  is not connected.
Define 
\[R_z \defeq \max \{r>0: (B(z,r) \cap {D^{\mathcal Q} }) \subset (D^{\mathcal Q})^z\}\,.\]
 Then \(R_z>0\). 
Moreover, since \(B(z,r)\subset B(z,s)\) for \(0<r<s\) and \(\cup_{r>0}B(z,r)=D\),
 it follows that \(R_z<+\infty\).  For every \(r>R_z\) the Kobayashi
 ball \(B(z,r)\) intersects \(D^ {\mathcal Q}\) in another
 connected component different from \((D^{\mathcal Q})^z\).
 Therefore, there
 exists \(w\in \de B(z,R_z)\) such that \(w\in {D^{\mathcal Q}}\) and \(w\not\in
 (D^{\mathcal Q})^z\). Let \( (D^{\mathcal Q})^w\)
 be the connected component of $w$ in \({D^{\mathcal Q}}\). As before, there 
exists a maximal \(R_w>0\) such that \(B(w,r)\cap {D^{\mathcal Q}} \subset (D^ {\mathcal Q})^w\).
Let \(A \defeq B(z,R_z)\cap B(w,R_w)\). By construction, \(A\cap {D^{\mathcal Q}}=\emptyset\). Moreover, since both 
\(B(z,R_z)\) and \(B(w,R_w)\) are \(\mathcal Q\)-invariant,
it follows that \(A\) is \(\mathcal Q\)-invariant.

Since \(A\) is a nonempty bounded convex domain by Lemma \ref{lemma:retraction}
and  \(\mathcal Q\) generates a relatively compact subgroup of \(\Aut A\),
there exists a point of \(A\) fixed by \(\mathcal Q\)
 (see \cite{Aba},  Thm. 2.5.7). Hence \(A\cap {D^{\mathcal Q}} \neq \emptyset\), giving a contradiction.
 \end{proof}
 
 \smallskip
 \noindent
We conclude with a result on the fixed point sets of the isotropy subgroups of $G$.

\smallskip
\begin{proposition}\label{rigidity}
Let  \(M \cong G/K\) be a hyperbolic homogeneous complex manifold and 
\(\mathcal F\) a $G$-invariant holomorphic foliation on \(M\).
Consider the fixed point set
\[
M^K=\{q\in M: k \cdot q=q \quad \hbox{for all } k \in K\}
\]
and the foliation $\mathcal F^K$ on $M^K$ induced by $\mathcal F$. Then
 the hyperbolic complex submanifold $M^K$ is connected and homogeneous with respect to a free  action of a Lie group $N$.
 Moreover, $\mathcal F^K$ is $N$-invariant 
and  every leaf of $\mathcal F^K$ 
 is given by  $F\cap M^K$, with $F$ a leaf of $\mathcal F$.
\end{proposition} 

\begin{proof}
By  \cite{Nak},  $M$ is biholomorphic to a convex domain in $\C^n$. Hence, by Proposition \ref{lemma:connected}, the fixed point set $M^K$ is a connected hyperbolic  complex manifold. Moreover, \[M^K= N_G(K) \cdot p= N \cdot p\,,\]
 where \(N \defeq N_G(K)^e/K \) and \(N_G(K)^e\) is the identity component in \(N_G(K)\), the normalizer of \(K\) in \(G\). 
 The leaves of \(\mathcal F^K\) are the 
 connected components of  \(F \cap M^K\), with \(F\) varying among the leaves of \(\mathcal F\).
 By construction, \(\mathcal F^K\) is 
 \(N\)-invariant. 
The fixed point set \(F^K\)  of \(K\) in \(F\) is connected by Proposition \ref{lemma:connected}.
Since \(F \cap M^K=F^K\), it follows that the intersection
\(F \cap M^K\) is a leaf of \(\mathcal F^K\).
\end{proof}

\appendix 
\section{} 
\label{Lie}

Here we give a different  proof of Proposition \ref{CLOSED} which exploits the
Lie group structure of the automorphism group of a Siegel domain.
The proof follows at once from the following lemmas.

\begin{lemma}
\label{SOLVABLESUB}
Let \(M\) be a hyperbolic homogeneous manifold and let \(G\) be  a closed subgroup of \(\Aut{M}\) which acts transitively on \(M\). 
Then \(G\) contains a closed,  simply connected, solvable subgroup \(S \) which acts
freely and transitively on \(M\).
\end{lemma}

\begin{proof}
Without loss of generality we may assume that \(G\) is connected. 
Thus it is contained in the identity component 
of the automorphism group of  \(M\). By \cite{Nak}, the 
manifold $M$
is biholomorphic to a homogeneous Siegel domain of type II.
By  \cite{Kan}, p. 38, Thm. B, the group $G$ admits a faithful finite dimensional representation.
Then from the Levi decomposition it follows that $G$ decomposes as a semidirect product
 $$G = (TL) \ltimes P\,,$$
 where \(P\) is a simply connected solvable Lie group, \(T\) is a compact torus and  \(L\) is a closed real semisimple 
 subgroup of $G$ centralizing \(T\) (see, e.g. \cite{Var} Ex. 44e, p.256).
Let \(KAN\) be an Iwasawa decomposition of \(L\) (cf. \cite{Hel}). 
We claim that the closed, simply connected, solvable subgroup $S:= (AN) \ltimes P$ 
of $G$ acts freely and transitively on $M$.

Indeed $P$ is contractible (\cite{Var}, Thm. 3.18.11), therefore so is $S$.
As a consequence, the group  $S$ admits no non-trivial compact subgroups and,  
since the $S$-action is proper, it  acts freely on $M$. 

Finally, $TK$ is a compact subgroup of $G$, thus it 
is contained in (in fact, it coincides with)
 an isotropy subgroup of the transitive Lie group $G$ (cf.  \cite{Kan}, Prop. 5.7).
As a consequence, $S$ acts 
transitively on $M$.
\end{proof}

\begin{remark} The above solvable subgroup $S$ is not necessarily split.
One can construct examples of  simply connected,  real non-split, solvable Lie groups 
acting  freely and transitively on certain Hermitian symmetric spaces.
\end{remark}  

\begin{lemma}
\label{SOLVABLE}
Let \(M\) be a hyperbolic homogeneous manifold 
endowed with a holomorphic foliation which is invariant with respect
to a closed, connected subgroup $S$ of $\Aut M$ acting freely and transitively on $M$.
Then
\begin{itemize}
\item[(i)] \(S\) is solvable and simply connected. In particular, every 
connected Lie subgroup of \(S\) is closed.
\item[(ii)] The leaves of the foliation are closed.
\end{itemize}
\end{lemma}

\begin{proof}
(i) From Lemma \ref{SOLVABLESUB}
it follows that $S$ is simply connected and solvable.
Thus every connected Lie subgroup of \(S\) is closed (cf. \cite{Var}, Thm. 3.18.12).

(ii) Let \(L_p\) be the stabilizer in $S$ of a leaf \(F_p\).
Then \(F_p\) is the orbit
of  \(L_p\).
Since \(M\) is diffeomorphic to \(S\)  via the $S$-action and
\(F_p\) is connected, it follows that \(L_p\) is connected.
Thus $L_p$ is closed in $S$ by 
(i) and  \(F_p\) is closed in $M$.
\end{proof}

%------------------------Bibliography-----------------------------


\begin{thebibliography}{BaBo}
\bibitem[Aba]{Aba} M.~Abate: {Iteration theory of holomorphic maps on taut manifolds}. {Research and Lecture Notes in Mathematics. Complex Analysis and Geometry}, {Mediterranean Press, Rende}, {1989}.
\bibitem[BaBo]{BaBo} P.~Baum, R.~Bott: {\sl Singularities of holomorphic foliations}. J. Differential Geometry 7 (1972) 279--342.
\bibitem[BeSc]{BeSc} A.~Behague, B.~Sc\'ardua: {\sl Foliations invariant under Lie group transverse actions}. Monatsh. Math. 153, 4 (2008) 295--308.
\bibitem[Boc]{Boc} S.~Bochner: {\sl Compact groups of differentiable transformations}. Annals of Math. 46, 3 (1945) 372-381.
\bibitem[BrSa]{BrSa} F.~Bracci, A.~Saracco: {\sl Hyperbolicity in unbounded convex domains}.   Forum Math. 5, 21 (2009) 815-826.
\bibitem[Bru]{Bru} M.~Brunella: {\sl Birational geometry of foliations}. Publicacoes Matem\'aticas do IMPA. [IMPA Mathematical Publications] Instituto de Matem\'atica Pura e Aplicada (IMPA), Rio de Janeiro, 2004.
%\bibitem[Che]{Che} S. S.~Chern: {\sl Complex manifolds without potential theory}. Second ed., %Universitext, Springer-Verlag, New York, 1995.
\bibitem[For]{For} F.~Forstneri\v{c}:  {\sl Interpolation by holomorphic automorphisms and embeddings in \(\C^n\)}. J. Geom. Anal. 9, 1 (1999), 93--117.
\bibitem[Gea]{Gea} L.~Geatti:
{\sl Holomorphic automorphisms of the tube domain over the Vinberg cone}.
Atti della Accademia Nazionale dei Lincei, Rendiconti, vol. LXXX (5), (1988) 283 -- 291.
\bibitem[Hel]{Hel}S.~Helgason:
{\sl Differential geometry, Lie groups and symmetric spaces.} 
 GSM {\bf 34}, AMS, Providence, 2001.
\bibitem[Ish]{Ish} H.~Ishi:
{\sl A torus subgroup of the isotropy group of a bounded homogeneous domain}.
Manuscripta Math. 130 (2009) 353--358.
\bibitem[Kan]{Kan} S.~Kaneyuki:
 {\sl Homogenous bounded domains
and Siegel domains}. Lect. Notes in Math. 241, Springer-Verlag, Berlin Heidelberg New York, 1971.
%\bibitem[Knapp]{Knapp} A.~W.~Knapp: 
%{\sl Representation theory of semisimple groups},
%{Princeton Landmarks in Mathematics},
%{Princeton University Press, Princeton, NJ}, {2001}.
\bibitem[Kob]{Kob} S.~Kobayashi:
{\sl Hyperbolic complex spaces}. Springer-Verlag, Berlin Heidelberg, 1998.
\bibitem[Mie]{Mie} C.~Miebach:
{\sl Quotients of bounded homogeneous domains by cyclic groups}. Osaka J. Math., 47 (2010) 331--352.
\bibitem[Nag]{Nag} S.~Nag:
{\sl Hyperbolic manifolds admitting holomorphic fiberings}. Bull. Austral. Math. Soc. 26, 2 (1982) 181--184.
\bibitem[Nak]{Nak} K.~Nakajima: {\sl Homogeneous hyperbolic manifolds and homogeneous Siegel domains}.
J. Math. Kyoto Univ. 25, 2 (1985) 269--291.
\bibitem[Roy]{Roy} H.~L.~Royden: {\sl Holomorphic fiber bundles with hyperbolic fiber}. Proceedings of the AMS. 42, 2 (1974) 311--312.
\bibitem[Tis]{Tis} D.~Tischler: {\sl On fibering certain foliated manifolds over \(S^1\)}. Topology 9 (1970) 153--154.
\bibitem[Var]{Var} V.~S.~Varadarajan:
{\sl Lie Groups, Lie Algebras, and Their Representation}.
GTM 102,  Springer-Verlag, Berlin, 1984.
%\bibitem[Vin1]{Vin1} E.~B.~Vinberg: {\sl Lie groups and Lie algebras, III}, 
%{Encyclopaedia of Mathematical Sciences},
%{41}, {Springer-Verlag, Berlin}, {1994}.
%\bibitem[Vin2]{Vin2} E.~B.~Vinberg:
%{\sl The Morozov-Borel theorem for real Lie groups.}
%Dokl. Akad. Nauk SSSR 141, 270-273 (1961)

%\bibitem[Xu]{Xu} Y. Xu, {\sl Theory of complex homogeneous bounded domains}, {Mathematics and its Applications}, {569}, {Translated and revised from the 2000 Chinese original},
 %{Science Press, Beijing; Kluwer Academic Publishers, Dordrecht}, {2005}.

\end{thebibliography}
\end{document}